\documentclass{amsart}

\usepackage[all]{xypic}
\usepackage[centertags]{amsmath}
\usepackage{amsfonts}
\usepackage{amscd}
\usepackage{amssymb}
\usepackage{amsthm}
\usepackage{newlfont}
\usepackage{amsxtra}
 \newtheorem{thm}{Theorem}[section]

 \newtheorem{lem}[thm]{Lemma}
 \newtheorem{prop}[thm]{Proposition}
 \theoremstyle{definition}
 \newtheorem{defn}[thm]{Definition}
 \theoremstyle{remark}
 \newtheorem{rem}[thm]{Remark}
 \theoremstyle{remark}
 
 \theoremstyle{definition}
 \newtheorem{notn}[thm]{Notation}
 \numberwithin{equation}{section}

 \newcommand{\Spf}{\mathrm{Spf}}

 \newcommand{\Aut}{\mathrm{Aut}}
 
 \newcommand{\Pic}{\mathrm{Pic}}
 
 \newcommand{\ord}{\mathrm{ord}}

 \newcommand{\Gal}{\mathrm{Gal}}
 \newcommand{\GL}{\mathrm{GL}}
 \newcommand{\PGL}{\mathrm{PGL}}
 \newcommand{\Fx}{\mathrm{Fix}}
 
 \newcommand{\sep}{\mathrm{sep}}

 \newcommand{\Fix}{\mathrm{Fix}}
 
 \newcommand{\Tr}{\mathrm{Tr}}

 \renewcommand{\mod}{\mathrm{mod}}
 \newcommand{\Nr}{\mathrm{Nr}}
 \newcommand{\fa}{\mathfrak a}
 
 \newcommand{\fb}{\mathfrak b}
 \newcommand{\fc}{\mathfrak c}

 \newcommand{\fp}{\mathfrak p}
 \newcommand{\fr}{\mathfrak r}
 
 \newcommand{\fn}{\mathfrak n}

 \newcommand{\fR}{\mathfrak R}
 
 \newcommand{\cO}{\mathcal{O}}

 \renewcommand{\cD}{\mathcal{D}}

 \newcommand{\C}{\mathbb{C}}
 
 \newcommand{\F}{\mathbb{F}}
 \newcommand{\M}{\mathbb{M}}
 \newcommand{\Q}{\mathbb{Q}}
 \newcommand{\Z}{\mathbb{Z}}

 \renewcommand{\P}{\mathbb{P}}

 \newcommand{\Fi}{F_\infty}

 \newcommand{\bg}{\overline{\Gamma}}
 \newcommand{\G}{\Gamma}
 
 \newcommand{\La}{\Lambda}

\begin{document}

\title{On hyperelliptic modular curves over function fields}

\author{Mihran Papikian}

\address{Department of Mathematics, Pennsylvania State University, University Park, PA 16802}

\email{papikian@math.psu.edu}

\thanks{The author was supported in part by NSF grant DMS-0801208 and Humboldt Research Fellowship.}



\begin{abstract}
We prove that there are only finitely many modular curves of
$\cD$-elliptic sheaves over $\F_q(T)$ which are hyperelliptic. In
odd characteristic we give a complete classification of such curves.
\end{abstract}


\maketitle


\section{Introduction}

The hyperelliptic classical modular curves and the hyperelliptic
Shimura curves over $\Q$ have been classified and extensively
studied, cf. \cite{Michon}, \cite{Ogg}, \cite{Ogg2}. One of the
reasons for the interest in hyperelliptic modular curves is that the
existence of a degree-$2$ map to the projective line makes such
curves amenable for explicit calculations.

Let $\F_q$ be the finite field with $q$ elements, $A=\F_q[T]$ be the
ring of polynomials in $T$ with $\F_q$-coefficients, and $F=\F_q(T)$
be the fraction field of $A$. To each ideal $\fn\lhd A$ there is an
associated Drinfeld modular curve $X_0(\fn)$ over $F$, which
classifies rank-$2$ Drinfeld $A$-modules with a cyclic $\fn$-torsion
subgroup. These curves are the analogues of classical modular curves
$X_0(N)$. Schweizer \cite{Schweizer} determined all $\mathfrak{n}$
for which $X_0(\mathfrak{n})$ is hyperelliptic, thus proving the
analogue of a result of Ogg \cite{Ogg}. The function field analogue
of Shimura curves was introduced by Laumon, Rapoport and Stuhler in
\cite{LRS}. These curves are moduli spaces of certain objects,
called $\cD$-elliptic sheaves, which generalize the notion of
Drinfeld modules. The purpose of the present paper is to determine
those modular curves of $\cD$-elliptic sheaves which are
hyperelliptic (Theorem \ref{thm-main}). This is the analogue of a
result of Michon and Ogg \cite{Michon}, \cite{Ogg2}.

\section{Notation}\label{Sec1}

Let $C:=\P^1_{\F_q}$ be the projective line over $\F_q$. Denote by
$F=\F_q(T)$ the field of rational functions on $C$. The set of
closed points on $C$ (equivalently, places of $F$) is denoted by
$|C|$. For each $x\in |C|$, we denote by $\cO_x$ and $F_x$ the
completions of $\cO_{C,x}$ and $F$ at $x$, respectively. The residue
field of $\cO_x$ is denoted by $\F_x$, the cardinality of $\F_x$ is
denoted by $q_x$, the degree $m$ extension of $\F_x$ is denoted by
$\F_x^{(m)}$, and $\deg(x):=\dim_{\F_q} (\F_x)$.

\vspace{0.1in}

Let $A:=\F_q[T]$ be the ring of polynomials in $T$ with $\F_q$
coefficients; this is the subring of $F$ consisting of functions
which are regular away from $\infty:=1/T$. The completion of an
algebraic closure of $F_\infty$ is denoted $\C_\infty$. Fox each
$x\in |C|-\infty$, we denote by $\fp_x\lhd A$ the corresponding
prime ideal of $A$. For $a\in A$, let $\deg(a)$ denote the degree of
$a$ as a polynomial in $T$. An ideal $\fa\lhd A$ is necessarily
principal $\fa=(a)$ for some $a\in A$; we define
$\deg(\fa)=\deg(a)$. Note that this does not depend on the choice of
a generator of $\fa$, and moreover $\deg(\fa)=\dim_{\F_q}(A/\fa)$
and  $\deg(\fp_x)=\deg(x)$.

\vspace{0.1in}

Given a ring $H$, we denote by $H^\times$ the group of its units.

\vspace{0.1in}

Let $D$ be a quaternion algebra over $F$, i.e., a $4$-dimensional
central simple algebra over $F$. For $x\in |C|$, we let
$D_x:=D\otimes_F F_x$. We assume throughout the paper that $D$ is
split at $\infty$, i.e., $D_\infty\cong \M_2(F_\infty)$. (Here
$\M_2$ is the ring of $2\times 2$ matrices.) Let $R$ be the set of
places where $D$ is ramified. It is a well-known fact that the
cardinality of $R$ is even, and conversely, for any choice of a
finite set $R\subset |C|$ of even cardinality there is a unique
quaternion algebra ramified exactly at the places in $R$; see
\cite[p. 74]{Vigneras}. If $R\neq \emptyset$, then $D$ is a division
algebra; if $R=\emptyset$, then $D\cong \M_2(F)$. The ideal
$\fr:=\prod_{x\in R}\fp_x$ is the \textit{discriminant} of $D$. Let
$D^\times$ be the algebraic group over $F$ defined by
$D^\times(B)=(D\otimes_F B)^\times$ for any $F$-algebra $B$; this is
the multiplicative group of $D$. Let $\alpha\mapsto \alpha'$ denote
the canonical involution of $D$; thus
$(\alpha\beta)'=\beta'\alpha'$. The \textit{reduced trace} of
$\alpha$ is $\Tr(\alpha)=\alpha+\alpha'$ and the \textit{reduced
norm} of $\alpha$ is $\Nr(\alpha)=\alpha\alpha'$.

\section{Preliminaries}\label{SecPrelim}

From now on we assume that $D$ is fixed and is a division algebra.
In particular, $R\neq \emptyset$. Fix a maximal $A$-order $\La$ in
$D$. Denote $\G:=\La^\times$. Since $D$ is split at $\infty$, it
satisfies the so-called \textit{Eichler condition} relative to $A$,
cf. \cite[(34.3)]{Reiner}. This implies that, up to conjugation,
$\La$ is the unique maximal $A$-order in $D$, i.e., any other
maximal $A$-order in $D$ is of the form $\alpha \La\alpha^{-1}$ for
some $\alpha\in D^\times(F)$, cf. \cite[Cor. III. 5.7]{Vigneras}.

Let $\Omega$ denote the Drinfeld upper half-plane over $\Fi$. As a
set $\Omega=\C_\infty-\Fi$, but $\Omega$ also has a natural
structure of a rigid-analytic space, cf. \cite{vdPut}. The group
$\G$ can be considered as a subgroup of $\GL_2(\Fi)$ via the
embedding
$$
\G\subset D^\times(F)\subset D^\times(\Fi)\cong \GL_2(\Fi).
$$
Hence $\G$ acts on $\Omega$ via linear fractional transformations.
Denote the image of $\G$ in $\PGL_2(\Fi)$ by $\bg$. $\bg$ is a
finitely generated, discrete and cocompact subgroup of
$\PGL_2(\Fi)$, so the quotient $\Omega/\G$ is the underlying
rigid-analytic space of a smooth projective geometrically
irreducible curve $X^R$ over $\Fi$; see \cite[Thm. 3.3]{vdPut}. In
fact, $X^R$ has a canonical model over $F$, since it is a coarse
moduli scheme of $\cD$-elliptic sheaves over $F$ with pole $\infty$,
where $\cD$ is a maximal $\cO_C$-order in $D$, cf. \cite[Ch. 4]{BS}
and \cite{PapGenus}. Using this moduli-theoretic interpretation, one
can show that the curve $X^R$ has good reduction at every $o\in
|C|-R-\infty$; see \cite[Cor. 3.2]{PapGenus}. We will denote the
genus of $X^R$ by $g(X^R)$ and the reduction of $X^R$ at $o$ by
$X^R_o$.

For a non-empty finite subset $S\subset |C|$, let
$$
\wp(S)=\left\{
         \begin{array}{ll}
           0, & \hbox{if some place in $S$ has even degree;} \\
           1, & \hbox{otherwise.}
         \end{array}
       \right.
$$

\begin{thm}\label{thm-genus}
$$
g(X^R)=1+\frac{1}{q^2-1}\prod_{x\in R}(q_x-1)-\frac{q}{q+1}\cdot
2^{\# R-1}\cdot \wp(R).
$$
\end{thm}
\begin{proof}
\cite[Thm. 5.4]{PapGenus}.
\end{proof}

\begin{thm}\label{thm-ss} For $o\in |C|-R-\infty$, we have
$$
\# X^R_o(\F_o^{(2)})\geq \frac{1}{q^2-1}\prod_{x\in R\cup
o}(q_x-1)+\frac{q}{q+1}\cdot 2^{\# R}\cdot \wp(R\cup o).
$$
\end{thm}
\begin{proof}
\cite[Cor. 4.8]{PapGenus}.
\end{proof}

Let $X$ be a curve over a field $K$. (From now on a curve is always
assumed to be smooth, projective and geometrically irreducible.)
Denote an algebraic closure of $K$ by $\bar{K}$, and let $\Aut(X)$
denote the group of $\bar{K}$-automorphisms of $X$. Denote
$X_{\bar{K}}:=X\otimes_{K}\bar{K}$.

\begin{thm}\label{thm-LKH} Let $X$ be a curve over $K$ of genus greater or equal to $2$.
The following conditions are equivalent:
\begin{enumerate}
\item There exists a $K$-morphism $w: X\to X'$ of
degree $2$, where $X'$ is a curve of genus $0$.
\item There exists a $\bar{K}$-morphism $w:X_{\bar{K}}\to
\P^1_{\bar{K}}$ of degree $2$.
\item $\Aut(X)$ contains an involution $\sigma$ defined over $K$ such that
the quotient $X/\sigma$ has genus $0$.
\item $\Aut(X)$ contains an involution $\sigma$ such that the
quotient $X_{\bar{K}}/\sigma$ is isomorphic to $\P^1_{\bar{K}}$.
\end{enumerate}
Moreover, the involution $\sigma$ is uniquely determined by $(4)$.
It is defined over $K$, and is in the center of $\Aut(X)$.
\end{thm}
\begin{proof}
See \cite[$\S$5]{LK}.
\end{proof}

\begin{defn} A curve $X$ satisfying the conditions of Theorem \ref{thm-LKH} is
called \textit{hyperelliptic}, and the involution $\sigma$ is called
the \textit{canonical involution} of $X$.
\end{defn}

\section{Main result} \label{SecProof}

\begin{thm}\label{thm-main}\hfill
\begin{enumerate}
\item For a fixed $q$ there are only finitely many $R$ such that $X^R$ is
hyperelliptic.
\item If $q$ is odd, then $X^R$ is hyperelliptic if and only if
$R=\{x,y\}$ and $$\{\deg(x), \deg(y)\}=\{1,2\}.$$
\end{enumerate}
\end{thm}

Most of the key ideas in the proof of this theorem go back to Ogg
\cite{Ogg}, \cite{Ogg2}. There are three main tools used in the
proof. Two of them are the genus formula for $X^R$ (Theorem
\ref{thm-genus}) and the estimate for the number of rational points
on a reduction of $X^R$ (Theorem \ref{thm-ss}). The third one is the
study of certain involutions on $X^R$, which are the analogues of
Atkin-Lehner involutions on Shimura curves. These involutions will
be discussed later in this section.

When $q$ is a power of $2$, to classify those (finitely many) $R$
for which $X^R$ is hyperelliptic one can follow the same strategy as
for odd $q$. Nevertheless, the argument will be technically much
more complicated. In fact, almost all preliminary results for the
proof of Part (2) of Theorem \ref{thm-main} crucially use the
assumption that the characteristic is not $2$.

\vspace{0.1in}

\noindent\textit{Proof of part (1) of Theorem \ref{thm-main}.}
Suppose $X^R$ is hyperelliptic. Fix some $o\in |C|-R-\infty$, and
consider $X^R_o$. By \cite[Prop. 5.14]{LK}, $X^R_o$ is also a
hyperelliptic curve. Over a finite field a curve of genus $0$ has a
rational point, hence is isomorphic to the projective line. In
particular, there is a degree-$2$ morphism $X^R_o\to \P^1_{\F_o}$
over $\F_o$. This implies
$$
\#X^R_o(\F_o^{(2)})\leq 2\#\P^1_{\F_o}(\F_o^{(2)}) =2(q_o^2+1).
$$
Using Theorem \ref{thm-ss}, we get
\begin{equation}\label{eq-fh}
\prod_{x\in R\cup o}(q_x-1)\leq 2(q_o^2+1)(q^2-1).
\end{equation}
Choose $o$ of minimal possible degree. It is enough to show that
(\ref{eq-fh}) can hold only for finitely many $R$ with such a choice
of $o$. Any place in $R$ of degree larger or equal to $\deg(o)$ only
increases the left hand-side of the inequality, so we can assume
that $R=\{x\in |C|-\infty\ |\ \deg(x)< \deg(o)\}$. For such $R$, the
left hand-side of (\ref{eq-fh}) as a polynomial in $q$ has degree
which is exponential in $\deg(o)$, but the right hand-side has
degree which is linear in $\deg(o)$. If there are infinitely many
$R$, then $\deg(o)$ can be arbitrarily large, which leads to a
contradiction. \hfill$\square$

\vspace{0.1in}

From now on we assume that the characteristic of $F$ is odd.

\begin{defn}
The \textit{normalizer} of $\La$ (respectively, $\G$) in $D$ is
$$
N(\La):=\{g\in D^\times(F)\ |\ g\La g^{-1}=\La\},
$$
respectively, $N(\G):=\{g\in D^\times(F)\ |\ g\G g^{-1}=\G\}$.
Clearly, $N(\La)\subseteq N(\G)$.

For a place $x\in |C|-\infty$, we define the local versions
$$
N(\La_x):=\{g\in D^\times(F_x)\ |\ g\La_x g^{-1}=\La_x\},
$$
where $\La_x=\La\otimes_A \cO_x$, and $ N(\La_x^\times):=\{g\in
D^\times(F_x)\ |\ g\La_x^\times g^{-1}=\La_x^\times\}$.
\end{defn}

\begin{prop}\label{prop-AL}\hfill
\begin{enumerate}
\item $N(\La)=\{g\in D^\times(F)\ |\ g\in N(\La_x),\ x\in
|C|-\infty\}$.
\item $N(\G)=\{g\in D^\times(F)\ |\ g\in N(\La_x^\times),\ x\in
|C|-\infty\}$.
\item $N(\La)=N(\G)$.
\item If $0\neq g\in \La$ and $\Nr(g)$ divides $\fr$, then $g\in N(\La)$.
\item $N(\La)/F^\times \G  \cong (\Z/2\Z)^{\# R}$.
\item There exists a set of elements $\{\pi_x\in \La\}_{x\in R}$ with
$(\Nr(\pi_x))=\fp_x$. Any such set generates $N(\La)/F^\times \G$.
\end{enumerate}
\end{prop}
\begin{proof} From the local-global correspondence for orders,
$g\La g^{-1}= \La$ if and only if $g\La_x g^{-1}=(g\La g^{-1})_x=
\La_x$ for all $x\in |C|-\infty$. This proves (1). A similar
argument also proves (2).

If $x\not\in R$, then $\La_x$ is conjugate to $\M_2(\cO_x)$. Since
the normalizer of $\M_2(\cO_x)$ in $D_x\cong \M_2(F_x)$ is
$F_x^\times \GL_2(\cO_x)$, we conclude $N(\La_x)=F_x^\times
\La_x^\times$. Since $\GL_2(\cO_x)$ contains an $\cO_x$-basis of
$\M_2(\cO_x)$, if $g\in D_x$ normalizes $\GL_2(\cO_x)$, then it also
normalizes $\M_2(\cO_x)$. Hence $N(\La_x^\times)=N(\La_x)=F_x^\times
\La_x^\times$. Next, if $x\in R$, then the uniqueness of the maximal
$\cO_x$-order in $D_x$ for $x\in R$ implies
$N(\La_x)=D^\times(F_x)$. Since $N(\La_x)\subseteq
N(\La_x^\times)\subseteq D^\times(F_x)$,
$N(\La_x^\times)=D^\times(F_x)$. Overall, we conclude that
$N(\La_x)=N(\La_x^\times)$ for all $x\in |C|-\infty$. Therefore, (3)
follows from (1) and (2).

Suppose $g\in \La$ and $(\Nr(g))$ divides $\fr$. Then $g\in
\La_x^\times$ for all $x\in |C|-R-\infty$, in particular, for such
$x$, $g\in N(\La_x)$. On the other hand, since
$N(\La_x)=D^\times(F_x)$ for $x\in R$, $g\in N(\La_x)$ for $x\in R$.
Overall, $g\in N(\La_x)$ for all $x\in |C|-\infty$, and (4) follows
from (1).

Let $x\in |C|-\infty$. Consider $N(\La_x)/F_x^\times \La_x^\times$.
As we discussed above, this quotient is trivial if $x\not\in R$. On
the other hand, when $x\in R$, the composition $\ord_x\circ \Nr$
induces an isomorphism $D^\times(F_x)/F_x^\times \La_x^\times\cong
\Z/2\Z$. Hence, for $x\in R$, $N(\La_x)/F_x^\times \La_x^\times\cong
\Z/2\Z$. By considering $g\in N(\La)$ as an element of $N(\La_x)$,
$x\in |C|-\infty$, we obtain a natural homomorphism
$$
N(\La)/F^\times \G \to \prod_{x\in |C|-\infty} N(\La_x)/F_x^\times
\La_x^\times\cong \prod_{x\in R}(\Z/2\Z).
$$
By (37.25) and (37.32) in \cite{Reiner}, this homomorphism is
surjective and has trivial kernel since $\Pic(A)=1$. This proves
(5).

The existence of the elements $\pi_x$ in (6) follows from
\cite[(34.8)]{Reiner}. The second claim of (6) follows from the fact
that the isomorphism in (5) is induced by $g \mapsto \prod_{x\in
R}\ord_x\circ \Nr(g)$.
\end{proof}

\begin{defn} Let $\fa\lhd A$ be an ideal dividing $\fr$, and let $\mu\in
\La$ be such that $(\Nr(\mu))=\fa$. By Proposition \ref{prop-AL},
such an element exists and $\mu\in N(\G)$. We can consider $\mu$ as
an element of $D^\times(\Fi)$, so $\mu$ acts on $\Omega$. Due to the
fact that it normalizes $\G$, $\mu$ induces an automorphism $w_\fa$
of $X^R$, cf. \cite[VII.1]{GvdP}. The ideal $(\mu)$ is a two-sided
integral ideal of $\La$ which is uniquely characterized by the fact
that $\Nr((\mu))=\fa$, cf. \cite[p. 86]{Vigneras}. Hence $w_\fa$
depends only on $\fa$, not on a particular choice of $\mu$.
Moreover, since $\mu^2\in F^\times \G$, $w_\fa^2=1$. We call $w_\fa$
the \textit{Atkin-Lehner involution} associated to $\fa$. By
Proposition \ref{prop-AL}, the Atkin-Lehner involutions form a
subgroup $W=N(\G)/F^\times \G$ of $\Aut(X^R)$ isomorphic to
$(\Z/2\Z)^{\# R}$. If $\fa$ and $\fb$ are divisors of $\fr$, and
$\fc=\mathrm{g.c.d.}(\fa,\fb)$, then $w_\fa w_\fb=w_{\fa\fb/\fc^2}$.
\end{defn}

The next theorem is the analogue of \cite[Thm. 2]{Rotger}.

\begin{thm}\label{thm-autAL}
If $R$ contains a place of even degree, then $\Aut(X^R)=W$.
\end{thm}
\begin{proof} Let $\widetilde{F}\neq F$ be a finite field extension of $F$.
From the theory of quaternion algebras it is well-known that
$\widetilde{F}$ embeds into $D$ as an $F$-subalgebra if and only if
$[\widetilde{F}:F]=2$ and no place in $R$ splits in $\widetilde{F}$.

Let $\gamma\in \G$ be a torsion element of order $n$: $\gamma^n=1$.
Let $p$ be the characteristic of $F$. If $p|n$, then
$\gamma^{n/p}\neq 1$ but $(\gamma^{n/p}-1)^p=0$. This is not
possible, since $D$ is a division algebra. We conclude that $\gamma$
is necessarily algebraic over $\F_q$. The field $\F_q(\gamma)\cong
\F_q^{(m)}$ is a finite extension of $\F_q$. Consider the subfield
$F(\gamma)$ of $D$. Since $[\F_q(\gamma):\F_q]=[F(\gamma):F]$, we
must have $m=1$ or $2$. We conclude that either $\gamma\in
\F_q^\times$, or $\F_{q^2}F$ embeds into $D$. A place $o\in |C|$ of
even degree splits in $\F_{q^2}F$. Therefore, if $R$ contains a
place of even degree then necessarily $\gamma\in \F_q^\times$. This
implies that the image $\bg$ of $\G$ in $\PGL_2(\Fi)$, which is a
finitely generated and discrete subgroup, is also torsion-free.
Hence $\bg$ is the Schottky group of $X^R$; see \cite{vdPut}. Let
$N_{\PGL_2(\Fi)}(\bg)$ denote the normalizer of $\bg$ in
$\PGL_2(\Fi)$. From the theory of Mumford curves \cite[p.
216]{GvdP}, one knows that $\Aut(X^R)\cong
N_{\PGL_2(\Fi)}(\bg)/\bg$.

Next, we claim that the $F$-vector space spanned by $\G$ is $D$. If
$\gamma\in \G$ is a non-torsion element, then $F(\gamma)$ is a
quadratic extension of $F$. Suppose $\gamma_1, \gamma_2\in \G$ are
such that $F(\gamma_1)\neq F(\gamma_2)$. Then one easily checks that
$\{1, \gamma_1, \gamma_2, \gamma_1\gamma_2\}$ are linearly
independent over $F$, hence form a basis of $D$. If such $\gamma_1$
and $\gamma_2$ do not exist, then $\G$ is an abelian group of rank
at most $1$. Such a group cannot have a cocompact image in
$\PGL_2(\Fi)$, which is a contradiction.

The previous paragraph implies that if $g\in N_{\GL_2(\Fi)}(\G)$,
then $g$ will actually normalize $D^\times(F)$. By the
Skolem-Noether theorem, $g$ induces an inner automorphism of $D$ so
that $g\in F_\infty^\times D^\times(F)$. Therefore,
$N_{\GL_2(\Fi)}(\G)=F_\infty^\times N(\G)$ and
$$\Aut(X^R)\cong N_{\PGL_2(\Fi)}(\bg)/\bg\cong N(\G)/F^\times\G=W.$$
\end{proof}

\begin{lem}\label{lem-AutPGL} Suppose $K$ is a field of characteristic
not equal to $2$. If $X$ is hyperelliptic, then $\Aut(X)$ does not
contain a subgroup isomorphic to $(\Z/2\Z)^4$.
\end{lem}
\begin{proof} We can assume $K$ is algebraically closed. The
canonical involution $\sigma$ is in the center of $\Aut(X)$, so if
$(\Z/2\Z)^4\subset \Aut(X)$, then $$(\Z/2\Z)^3 \subset
\Aut(X/\sigma)\cong \Aut(\P^1_K)\cong \PGL_2(K).$$ But $\PGL_2(K)$
does not contain a subgroup isomorphic to $(\Z/2\Z)^3$ if the
characteristic is not $2$, cf. \cite[p. 33]{Ledet}.
\end{proof}

\begin{prop}\label{prop-HEP}
If $X^R$ is hyperelliptic, then $R=\{x,y\}$ and $\{\deg(x),
\deg(y)\}$ is one of the following pairs
\begin{itemize}
\item $\{1,2\}, \{1,3\}, \{2,2\}$ if $q=3$;
\item $\{1,2\}$ if $q\geq 5$.
\end{itemize}
\end{prop}
\begin{proof} We know that $$W\cong (\Z/2\Z)^{\# R}\subset \Aut(X^R),$$ so if
$X^R$ is hyperelliptic, then Lemma \ref{lem-AutPGL} forces $\# R=2$.
Hence we can assume $R=\{x,y\}$ from now on. Now we apply the
argument in the proof of Part (1) of Theorem \ref{thm-main}. Note
that we can choose $o$ of degree $1$, since the number of degree $1$
places in $|C|-\infty$ is $q>2$. Then the following inequality must
hold:
$$
(q_x-1)(q_y-1)+4q\cdot \wp(R)\leq 2(q^2+1)(q+1).
$$
One easily checks that this is possible if and only if either
$$\{\deg(x),\deg(y)\}=\{1,1\}, \{1,2\}$$ and $q$ is arbitrary, or
$q=3$ and $$\{\deg(x),\deg(y)\}=\{1,1\}, \{1,2\}, \{1,3\},
\{2,2\}.$$ In $\{1,1\}$ case, $X^R$ has genus zero by Theorem
\ref{thm-genus}, so it is not hyperelliptic.
\end{proof}

To determine which of the curves listed in Proposition
\ref{prop-HEP} are actually hyperelliptic, we need to compute the
number of fixed points of Atkin-Lehner involutions acting on $X^R$.

\begin{defn}[cf. \cite{Vigneras}]
Let $L$ be a quadratic extension of $F$, and let $\fR$ be an
$A$-order in $L$ (e.g., the integral closure of $A$ in $L$). Suppose
$\varphi:L\hookrightarrow D$ is an embedding. One says that
$\varphi$ is an \textit{optimal embedding} of $\fR$ into $\La$ if
$\varphi(L)\cap \La=\varphi(\fR)$. Two optimal embeddings $\varphi$
and $\psi$  of $\fR$ into $\La$ are \textit{equivalent} if there is
$\gamma\in \G$ such that $\psi=\gamma \varphi \gamma^{-1}$. Denote
by $\Theta(\fR, \La)$ the number of inequivalent optimal embeddings
of $\fR$ into $\La$. By a result of Eichler \cite[pp.
92-94]{Vigneras}, we have
\begin{equation}\label{eq-Eich}
\Theta(\fR, \La)= h(\fR)\prod_{x\in
R}\left(1-\left(\frac{L}{x}\right)\right),
\end{equation}
where $h(\fR)$ is the class number of $\fR$ and
$\left(\frac{L}{x}\right)$ is the \textit{Artin-Legendre symbol}:
$$
\left(\frac{L}{x}\right)=\left\{
                           \begin{array}{ll}
                             1, & \hbox{if $x$ splits in $L$;} \\
                             -1, & \hbox{if $x$ is inert in $L$;} \\
                             0, & \hbox{if $x$ ramifies in $L$.}
                           \end{array}
                         \right.
$$
\end{defn}
\begin{defn}
A quadratic extension $L$ of $F$ is \textit{imaginary} if $\infty$
does not split in $L$.
\end{defn}

\begin{lem}
Let $f\in A$ be a polynomial of degree $d$ and let $c\in
\F_q^\times$ be its leading coefficient, i.e.,
$f=cT^d+c_{d-1}T^{d-1}+\cdots$. Then $L:=F(\sqrt{f})$ is imaginary
if and only if one of the following holds:
\begin{enumerate}
\item $d$ is odd;
\item $d$ is even and $c$ is not a square in $\F_q^\times$.
\end{enumerate}
In the first case $\infty$ ramifies in $L/F$, and in the second case
it remains inert. The integral closure of $A$ in $L$ is
$A[\sqrt{f}]$, which is also the subring of $L$ consisting of
functions regular away from $\widetilde{\infty}$ - the unique place
over $\infty$.
\end{lem}
\begin{proof}
See \cite[p. 187]{Lorenzini}.
\end{proof}

\begin{notn} Let $\fa|\fr$ and $f\in A$ be the monic generator of the ideal
$\fa$. Let $\kappa\in \F_q^\times$ be a fixed element which is not a
square in $\F_q^\times$. Denote $\fR_\fa:=A[\sqrt{\kappa f}]$ and
$\fR_\fa':=A[\sqrt{f}]$, and their fields of fractions by
$K_\fa=F(\sqrt{\kappa f})$ and $K_\fa':=F(\sqrt{f})$. Denote by
$\Fx(w_\fa)$ the set of fixed points of $w_\fa$ acting on $X^R$.
\end{notn}

\begin{prop}\label{prop-fp} With above notation,
$$
\#\Fx(w_\fa)=\left\{
               \begin{array}{ll}
                 \Theta(\fR_\fa, \La), & \hbox{if $\deg(f)$ is even;} \\
                 \Theta(\fR_\fa, \La)+\Theta(\fR_\fa', \La), & \hbox{if $\deg(f)$ is odd.}
               \end{array}
             \right.
$$
\end{prop}
\begin{proof} (Cf. \cite[$\S$2]{Ogg2}.)
We will show that $\Fx(w_\fa)$ is in one-to-one correspondence with
the inequivalent optimal embeddings of $\fR_\fa$ (and $\fR_\fa'$)
into $\La$ when $\deg(f)$ is even (resp. odd).

Let $\varphi$ be an optimal embedding of $\fR_\fa$ into $\La$. Then
$\mu:=\varphi(\sqrt{\kappa f})\in \La$. Note that $\mu'=-\mu$, since
the canonical involution of $D$ restricted to $F(\mu)$ induces the
non-trivial automorphism of this field over $F$. In particular,
$\Tr(\mu)=0$ and $\Nr(\mu)=-\kappa f$. Therefore, the action of
$\mu$ on $\Omega$ induces the action of $w_\fa$ on $X^R$. We claim
that $\mu$ has a fixed point in $\Omega$. Let $\begin{pmatrix} a & b
\\ c & d\end{pmatrix}$ be the matrix representation of $\mu$ under
the embedding $\La\hookrightarrow D(\Fi)\cong \M_2(\Fi)$. The fixed
points of $\mu$ correspond to the solutions of
\begin{equation}\label{eq-quadr}
c z^2 +(d-a)z -b =0.
\end{equation}
Since the matrix of $\mu$ has trace $0$ and determinant $-\kappa f$,
the discriminant of this quadratic equation is $4\kappa f$. Hence
its solutions generate an imaginary quadratic extension of $F$, in
particular, they are in $\C_\infty-\Fi$. The image of such a
solution on $X^R$ gives a fixed point of $w_\fa$. The same argument
also applies to an optimal embedding of $\fR_\fa'$ when $\deg(f)$ is
odd.

Now let $P\in X^R(\C_\infty)$ be a fixed point of $w_\fa$. This
means that there exists $z\in \Omega$ and $\mu\in \La$ with
$(\Nr(\mu))=\fa$, such that $\mu z = \gamma z$ for some $\gamma\in
\G$. By replacing $\mu$ by $\gamma^{-1}\mu$, we can assume $\mu
z=z$. Using a matrix representation of $\mu$ and (\ref{eq-quadr}),
this easily implies that the field extension $F(\mu)$ is necessarily
imaginary and $[\Fi(z):\Fi]=2$. Next, by Proposition \ref{prop-AL},
$\mu$ generates a two sided $\La$-ideal $I(\fa)=\mu\La=\La\mu$. By
\cite[p. 86]{Vigneras}, any $\tilde{\mu}$ with
$(\Nr(\tilde{\mu}))=\fa$ generates the same two-sided ideal
$I(\fa)$. In particular, $\mu'$ generates $I(\fa)$, so $\mu'=c\mu$
for some $c\in F(\mu)\cap \G=\F_q^\times$. (In fact, $c=-1$, since
$\mu+\mu'\in F$.) Hence $\mu^2=s f$, for some $s\in \F_q^\times$
such that $F(\sqrt{sf})$ is imaginary. Therefore, $A[\mu]$ is
isomorphic to $\fR_\fa$ (or $\fR_\fa'$ when $\deg(f)$ is odd), and
we obtain two optimal embedding of the corresponding ring into $\La$
given by $\mu$ and $\mu'$.

The existence of the rigid-analytic uniformization $\pi: \Omega\to
X^R(\C_\infty)$ is a consequence of a stronger result which shows
that the formal schemes $\hat{\Omega}/\G$ and $\hat{X}^R$ represent
the same functor over $\Spf(\cO_\infty)$, where $\hat{\Omega}$ is
Mumford's formal scheme with generic fibre $\Omega$ and $\hat{X}^R$
is the formal completion of $X^R$ along its closed fibre at
$\infty$; see \cite[Ch. 4]{BS}. This stronger statement implies that
$\pi$ is $\Gal(\Fi^\sep/\Fi)$-equivariant, in the sense that if
$z\in \Fi^\sep$ and $g\in \Gal(\Fi^\sep/\Fi)$, then $\pi(g
z)=g\pi(z)$.

So far to a fixed point $P$ of $w_\fa$ we have associated $z\in
\Omega$ in a quadratic extension of $\Fi$ and two optimal embeddings
of $\fR_\fa$ (or $\fR_\fa'$) into $\La$, corresponding to $\mu$ and
$\mu'$, and we showed that any optimal embedding arises in this
manner. The points in the orbit $\G z$ produce equivalent embeddings
since if $z$ corresponds to $\mu$ and $\mu'$, then $\gamma z$
corresponds to $\gamma \mu \gamma^{-1}$ and $\gamma \mu'
\gamma^{-1}$. Let $\tau$ be the generator of $\Gal(K_\fa/F)$. Note
that $\mu'=\mu\circ \tau$. Denote the equivalence class of an
optimal embedding $\mu$ by $[\mu]$. With this notation, to $P$ there
are two associated equivalence classes of embedding $\{[\mu],
[\mu\circ \tau]\}$. The auxiliary point $z$ in this construction
gives an embedding $K_\fa\hookrightarrow \C_\infty$, and allows us
to consider $\tau$ as an element of $\Gal(\Fi^\sep/\Fi)$. By the
$\Gal(\Fi^\sep/\Fi)$-equivariance of $\pi$, if $P=\pi(z)$ then $\tau
(P)=\pi(\tau z)$. Since $w_\fa$ is defined over $\Fi$, $\tau(P)$ is
also a fixed point of $w_\fa$. Note that the optimal embeddings
corresponding to $\tau(z)$ are $\mu\circ \tau$ and $\mu\circ
\tau\circ \tau=\mu$. Therefore, the equivalence classes of
embeddings corresponding to $\tau(P)$ are also $\{[\mu], [\mu\circ
\tau]\}$. Finally, it is easy to see that $[\mu]$ and $[\mu\circ
\tau]$ appear only among the embeddings associated to $P$ and
$\tau(P)$. Overall, we get a correspondence between $\Fix(w_\fa)$
and the equivalence classes of optimal embeddings of $\fR_\fa$ (and
$\fR_\fa'$) into $\La$:
$$
\{P, \tau P\}\leftrightarrow \{[\mu], [\mu\circ \tau]\}.
$$
To deduce the proposition it remains to show that $P=\tau(P)$ if and
only if $[\mu]=[\mu\circ \tau]$. Now, $P=\tau(P)$ if and only if
$\tau(z)=\gamma (z)$ for some $\gamma\in \G$. But $\tau(z)$ is a
fixed point of $\mu'$ and $\gamma (z)$ is a fixed point of
$\gamma\mu\gamma^{-1}$. Therefore, $P=\tau(P)$ if and only if
$\mu'=\gamma\mu\gamma^{-1}$, i.e., $[\mu\circ\tau]=[\mu]$.
\end{proof}

Let $a\in A$. The field $F(\sqrt{a})$ is the function field of the
curve $C'$ over $\F_q$ given by $y^2=a$. A genus formula for such a
curve is well-known (cf. \cite[p. 332]{Lorenzini}):
$$
g(C')=\left\lfloor \frac{\deg(a)-1}{2}\right\rfloor.
$$
Let $J_{C'}$ be the Jacobian variety of $C'$. The following fact is
also well-known (cf. \cite[Ch. VIII]{Lorenzini}):
\begin{lem}\label{lem-cn}
Suppose $F(\sqrt{a})$ is imaginary. Then
$$
h(A[\sqrt{a}])=\left\{
                 \begin{array}{ll}
                   \# J_{C'}(\F_q), & \hbox{if $\deg(a)$ is odd;} \\
                   2\# J_{C'}(\F_q), & \hbox{if $\deg(a)$ is even.}
                 \end{array}
               \right.
$$
\end{lem}

To simplify the notation, for $x, y\in R$ denote $w_x:=w_{\fp_x}$,
$w_y:=w_{\fp_y}$, $w_{xy}:=w_xw_y$. Let $f_x$ and $f_y$ be the monic
generators of $\fp_x$ and $\fp_y$, respectively.

\begin{lem}\label{lem1.16}
Suppose $R=\{x, y\}$, $\deg(x)=1$ and $\deg(y)=2$. Then
$$
\#\Fx(w_x)=0, \quad \#\Fx(w_y)=4, \quad \#\Fx(w_{xy})=2q+2.
$$
\end{lem}
\begin{proof}
From (\ref{eq-Eich}), Proposition \ref{prop-fp} and Lemma
\ref{lem-cn}, we have
$$
\#\Fx(w_x)=1\cdot
\left(1-\left(\frac{F(\sqrt{f_x})}{y}\right)\right)+1\cdot
\left(1-\left(\frac{F(\sqrt{\kappa f_x})}{y}\right)\right)=0.
$$
Similarly,
$$
\#\Fx(w_y)=2\cdot \left(1-\left(\frac{F(\sqrt{\kappa
f_y})}{x}\right)\right)=4.
$$
Let $f\in A$ be a polynomial of degree $3$. Then $E:z^2=f$ is an
elliptic curve over $\F_q$ (note that $\widetilde{\infty}$ is
rational on $E$). Since $E$ is its own Jacobian, $h(A[\sqrt{f}])$ is
equal to the number of solutions of $z^2=f$ over $\F_q$ plus $1$.
Now let $f=f_xf_y$. Then
\begin{align*}
\#\Fx(w_{xy}) &= h(A[\sqrt{f}]) + h(A[\sqrt{\kappa f}])\\
& =2
+\#\{(\alpha,\beta)\in \F_q\times \F_q\ |\ \alpha^2=f(\beta)\}\\
&+\#\{(\alpha,\beta)\in \F_q\times \F_q\ |\ \alpha^2=\kappa
f(\beta)\}.
\end{align*}
If $\beta\in \F_q$ is fixed and $f(\beta)\neq 0$, then exactly one
of $f(\beta)$ and $\kappa f(\beta)$ is a square in $\F_q$, and we
get two solutions of the corresponding equation. On the other hand,
if $f(\beta)=0$, then both $z^2=f(\beta)$ and $z^2=\kappa f(\beta)$
have exactly one solution, namely $z=0$. Overall, we get $2q$
solutions.
\end{proof}

\begin{lem}\label{lem1.17}
Suppose $R=\{x, y\}$, and $\deg(x)=\deg(y)=2$. Then
$$
\#\Fx(w_x)\leq 4, \quad \#\Fx(w_y)\leq 4, \quad \#\Fx(w_{xy})\leq
2(q+1+2\sqrt{q}).
$$
\end{lem}
\begin{proof}
The argument is similar to the proof of the previous lemma. Note
that $z^2=f_xf_y$ corresponds to a curve of genus $1$. The Jacobian
of such a curve is an elliptic curve, hence has at most
$(q+1+2\sqrt{q})$ rational points over $\F_q$ by the Hasse-Weil
bound.
\end{proof}

\begin{lem}\label{lem1.18}
Suppose $q=3$ and $R=\{x, y\}$, where $\deg(x)=1$ and $\deg(y)=3$.
Then
$$
\#\Fx(w_x)= 2, \quad \#\Fx(w_y)\geq 2, \quad \#\Fx(w_{xy})\geq 4.
$$
If $\#\Fx(w_y)= 2$, then $\#\Fx(w_{xy})=8$.
\end{lem}
\begin{proof}
We have
$$
\#\Fx(w_x)=1\cdot
\left(1-\left(\frac{F(\sqrt{f_x})}{y}\right)\right)+1\cdot
\left(1-\left(\frac{F(\sqrt{\kappa f_x})}{y}\right)\right).
$$
Since $[\F_y:\F_q]=3$, $\kappa$ is not a square in $\F_y$, so
exactly one of $f_x$ or $\kappa f_x$ is a square modulo $y$. This
implies $\#\Fx(w_x)=2$.

Now consider the curves defined by $z^2=f_y$, $z^2=\kappa f_y$ and
$z^2=\kappa f_xf_y$. Denote the Jacobians of these curves by $E_1$,
$E_2$ and $E_3$, respectively. All three of these abelian varieties
are elliptic curves. We have
$$
\#\Fx(w_y)=\#E_1(\F_q)\cdot
\left(1-\left(\frac{F(\sqrt{f_y})}{x}\right)\right)+
\#E_2(\F_q)\cdot \left(1-\left(\frac{F(\sqrt{\kappa
f_y})}{x}\right)\right).
$$
Since $\deg(x)=1$, $\kappa$ is not a square in $\F_x\cong \F_q$. On
the other hand, since $f_y$ is irreducible, $f_y\not \equiv 0\
(\mod\ x)$. Thus, exactly one of $f_y$ or $\kappa f_y$ is a square
modulo $x$, so $\#\Fx(w_y)\geq 2$. Suppose $\#\Fx(w_y)= 2$. This is
possible only if either $\#E_1(\F_q)=1$ or $\#E_2(\F_q)=1$. A
case-by-case verification shows that we must have
\begin{equation}\label{eq-f}
f_y=T^3-T+1 \quad \text{or}\quad f_y=T^3-T-1:
\end{equation}
if $f_y=T^3-T+1$ then $\#E_2(\F_q)=1$; if $f_y=T^3-T-1$ then
$\#E_1(\F_q)=1$. In both cases we indeed get $\#\Fx(w_y)=2$.

Next, $\#\Fx(w_{xy})=2\cdot \#E_3(\F_q)$. Note that $z^2=\kappa
f_xf_y$ has a rational solution: if $f_x=T-c$ with $c\in \F_q$, then
$T=c$, $z=0$ is a solution. Hence, $E_3$ besides its $0$ for the
group structure has another rational point, so $\#E_3(\F_q)\geq 2$.
This implies the desired lowed bound on $\#\Fx(w_{xy})$. Finally,
suppose $\#\Fx(w_y)=2$. From (\ref{eq-f}) we get that the function
$\kappa f_y(s)\in \F_3$, $s\in \F_3$, is either identically $1$ or
$-1$, and using this, one easily checks that $\#E_3(\F_q)=4$.
\end{proof}

\begin{lem}\label{lem1.19}
Suppose $X$ is a hyperelliptic curve over a field whose
characteristic is not $2$, and the genus of $X$ is even. Then any
involution of $X$, except for the canonical involution, has exactly
$2$ fixed points.
\end{lem}
\begin{proof}
If the characteristic is not $2$, then the proof of \cite[Prop.
1]{Ogg} applies.
\end{proof}

\vspace{0.1in}

\noindent\textit{Proof of part (2) of Theorem \ref{thm-main}.} Since
the characteristic of $F$ is odd, $X^R$ is hyperelliptic if and only
if $\Aut(X^R)$ contains an involution having $2g(X^R)+2$ fixed
points (by the Riemann-Hurwitz formula). If $R=\{x,y\}$ and
$\{\deg(x), \deg(y)\}=\{1,2\}$, then $g(X^R)=q$ (Theorem
\ref{thm-genus}). On the other hand, by Lemma \ref{lem1.16},
$\#\Fx(w_{xy})=2q+2$. Thus, $X^R$ is hyperelliptic in this case,
with $w_{xy}$ being the canonical involution. Thanks to Proposition
\ref{prop-HEP}, to conclude that this is the only possibility, we
need to check that when $q=3$, $R=\{x,y\}$ and $\{\deg(x),
\deg(y)\}=\{2,2\}, \{1,3\}$, $X^R$ is not hyperelliptic.

Suppose $q=3$ and $\{\deg(x), \deg(y)\}=\{2,2\}$. From Theorem
\ref{thm-autAL} we know that the canonical involution, if it exists,
is an Atkin-Lehner involution. Since $g(X^R)=9$, the canonical
involution must have $20$ fixed points. On the other hand, from
Lemma \ref{lem1.17} an Atkin-Lehner involution has at most $14$
fixed points. Hence $X^R$ is not hyperelliptic.

Now suppose $q=3$ and $\{\deg(x), \deg(y)\}=\{1,3\}$. Then
$g(X^R)=6$. Since $\#\Fx(w_{xy})\geq 4$, from Lemma \ref{lem1.18}
and Lemma \ref{lem1.19}, we deduce that if $X^R$ is hyperelliptic
then $w_{xy}$ is necessarily the canonical involution and
$\#\Fx(w_{y})=2$. If this is the case, then $\#\Fx(w_{xy})=8$. On
the other hand, the canonical involution must have $14$ fixed
points, a contradiction. \hfill $\square$

\vspace{0.1in}

\begin{rem}
As we mentioned in the introduction, the hyperelliptic Drinfeld
modular curves $X_0(\fn)$ over $F$ were classified by Schweizer
\cite{Schweizer}. When $q>2$ the answer is similar to Theorem
\ref{thm-main}: $X_0(\fn)$ is hyperelliptic if and only if
$\deg(\fn)=3$.
\end{rem}

\subsection*{Acknowledgments} I thank Professor E.-U. Gekeler for useful discussions.
The article was written while I was visiting Saarland University. I
thank the department of mathematics for its hospitality.



\end{document}